\documentclass[11pt,a4paper,reqno]{amsart}
\usepackage[margin=2.5cm]{geometry}
\usepackage{xcolor}

\allowdisplaybreaks

\usepackage[all,pdf]{xy}
\usepackage{amsmath}
\usepackage{amssymb}
\usepackage{amscd}
\usepackage{graphicx}
\usepackage{comment}
\usepackage{enumerate}

\usepackage[colorlinks=true]{hyperref}

\hypersetup{urlcolor=red, citecolor=blue}
\makeatletter
\@namedef{subjclassname@2020}{%
	\textup{2020} Mathematics Subject Classification}
\makeatother

\theoremstyle{plain}
\newtheorem{thm}{Theorem}[section]
        \newtheorem*{thm*}{Theorem}
\newtheorem{cor}[thm]{Corollary}
\newtheorem{lem}[thm]{Lemma}
\newtheorem{prop}[thm]{Proposition}

\newtheorem*{conj*}{Conjecture}

\theoremstyle{definition}

\newtheorem{rem}[thm]{Remark}

\def \Z {\mathbb Z}
\def \R {\mathbb R}

\def \M {\mathcal M}

\def\diam{{\mathrm{diam}}}
\def\epsilon{\varepsilon}

\numberwithin{equation}{section}

\newcommand{\Lip}{\mathrm{Lip}}

\title[Sharpness of the Dini condition]{Sharpness of the Dini condition in Liv\v{s}ic's theorem and the Ma\~n\'e lemma}


\date{}

\begin{document}
\author[O. Jenkinson and  H. Zhao]{Oliver Jenkinson and Hang Zhao}

\address[O.~Jenkinson]{School of Mathematical Sciences, Queen Mary University of London, Mile End Road,
London, E1 4NS, UK}
\email{o.jenkinson@qmul.ac.uk}

\address[H. Zhao]{School of Mathematics and Big Data, Neijiang Normal University,
Neijiang, Sichuan 641100, China}
\email{hangzhao@njtc.edu.cn, hangzhaoscu@163.com}
\thanks{The authors contributed equally to this work.}

\begin{abstract}
Liv\v{s}ic remarked that his theorem on the cohomological equation over hyperbolic systems remains valid, with a continuous solution, for functions whose modulus of continuity $\omega$ satisfies the Dini condition. We show that this condition cannot be weakened: if a concave modulus $\omega$ does not satisfy it, then on every infinite basic set of a $C^1$ diffeomorphism, some function with modulus of continuity a multiple of $\omega$ has zero integral for every invariant measure but does not have any bounded sub-action. So neither Liv\v{s}ic's theorem nor the Ma\~n\'e lemma extends to such moduli. The function is an explicit uniformly convergent series of coboundaries, defined for any Lipschitz map with an orbit whose past and future have disjoint closures. On such basic sets, the largest Birkhoff sums of length $N$ of functions with modulus of continuity $\omega$ and zero integrals are comparable to $\int_{e^{-N}}^1\omega(r)/r\,dr$, uniformly in $\omega$ and $N$.
\end{abstract}
\maketitle

\section{Introduction}\label{sec:intro}

The celebrated Liv\v{s}ic theorem \cite{Liv71,Liv72} states that, over a transitive hyperbolic system, a H\"older continuous function is a coboundary if and only if its sum over every periodic orbit vanishes. It has applications to rigidity (see e.g.~\cite{GK80,LMM86}), and has been extended in various directions, for instance to matrix cocycles \cite{Ka11} and to partially hyperbolic diffeomorphisms \cite{Wi13}.

 Liv\v{s}ic remarked that the H\"older condition can be weakened: it is enough that the modulus of continuity $\omega$ of the function satisfies the Dini condition $\int_0^1\omega(r)/r\,dr<\infty$. In this paper we show that the Dini condition cannot be weakened further.
 For every concave modulus of continuity $\omega$ which does not satisfy it, and every basic set of a $C^1$ diffeomorphism other than a periodic orbit, we write down a function with modulus of continuity a multiple of $\omega$ whose sum over every periodic orbit vanishes, but which is not a continuous coboundary.
The same function shows that the Dini condition is also sharp in the Ma\~n\'e lemma, the analogue of Liv\v{s}ic's theorem for cohomological inequalities and a basic tool in ergodic optimization.
We also show that, on these basic sets, the largest Birkhoff sums of length $N$ of functions with modulus of continuity $\omega$ whose sums over periodic orbits vanish are comparable to the truncated Dini integral $\int_{e^{-N}}^1\omega(r)/r\,dr$, uniformly in $\omega$ and $N$.

To state the results, let $T:X\to X$ be a continuous map of a non-empty compact metric space $(X,d)$, and let $\M_T$ be the set of $T$-invariant Borel probability measures. A \emph{continuous coboundary} is a function of the form $v\circ T-v$, where $v\in C(X)$ is called a \emph{transfer function}. The integral of a continuous coboundary vanishes for every $\mu\in\M_T$, and Liv\v{s}ic's theorem \cite{Liv71,Liv72} provides a converse for hyperbolic systems: if, for example, $T$ is a transitive Anosov diffeomorphism, and $f$ is a H\"older continuous function such that $\int f\,d\mu=0$ for every $\mu\in\M_T$ (equivalently, such that the sum of $f$ over each periodic orbit vanishes), then $f=v\circ T-v$ for some H\"older continuous function $v$ (see \cite[Theorem~19.2.1]{KH}).

In ergodic optimization (see e.g.~\cite{Boc,Je}), the role of the cohomological equation
$f=v\circ T -v$ is played by an inequality,
and that of Liv\v{s}ic's theorem by
the Ma\~n\'e lemma (also called the Ma\~n\'e--Conze--Guivarc'h lemma; see e.g.~\cite{Bo11, Je} for its history).
For $f\in C(X)$, let $S_nf:=\sum_{i=0}^{n-1}f\circ T^i$ for $n\ge1$, and
$\beta(f)=\beta(X,T,f):=\max_{\mu\in\M_T}\int f\,d\mu$.
A \emph{sub-action} for $f$ is a function $u:X\to\R$, not necessarily continuous, such that
\begin{equation}\label{eq:subaction}
f+u-u\circ T\le\beta(f).
\end{equation}
For hyperbolic systems, such as transitive Anosov diffeomorphisms, the Ma\~n\'e lemma states that every H\"older continuous function has a continuous sub-action \cite{Bo,LT}, and a Lipschitz one if the function is Lipschitz \cite{Bo11,STY}; for the case of expanding maps, see \cite{Bo00,CLT}.
 A measure $\mu\in\M_T$ is \emph{$f$-maximizing} if $\int f\,d\mu=\beta(f)$; if $u$ is a continuous sub-action for $f$, then $\mu\in\M_T$ is $f$-maximizing if and only if its support is contained in the closed set where equality holds in \eqref{eq:subaction}.

This paper concerns the regularity of $f$ needed in Liv\v{s}ic's theorem and the Ma\~n\'e lemma. We measure regularity by moduli of continuity: a \emph{concave modulus} is a continuous, non-decreasing, concave function $\omega:[0,\infty)\to[0,\infty)$ with $\omega(0)=0$ and $\omega(r)>0$ for $r>0$, and $\Lip_\omega(X)$ is the space of functions $f:X\to\R$ for which $|f|_\omega:=\sup_{x\ne y}|f(x)-f(y)|/\omega(d(x,y))$ is finite; for $\omega(r)=r^\alpha$, with $0<\alpha\le1$, this is the space of $\alpha$-H\"older functions. The modulus $\omega$ satisfies the \emph{Dini condition}, or is \emph{Dini}, if $\int_0^1\omega(r)/r \,dr<\infty$, and in that case its \emph{Dini transform} \cite{JZZ} is the concave modulus $D\omega(r):=\int_0^r\omega(s)/s \,ds$; for example, a concave modulus with $\omega(r)=(\log(1/r))^{-p}$ for all small $r>0$ is Dini if and only if $p>1$, and then $D\omega(r)=(\log(1/r))^{1-p}/(p-1)$ for all small $r>0$. 
Functions with a modulus of continuity satisfying the Dini condition are often said to have \emph{summable variation}, a term from symbolic dynamics (see Remark~\ref{rem:sumvar}); for the role of the Dini condition in thermodynamic formalism and in other questions of ergodic optimization, see e.g.\ \cite{FJ01,JZZ,JLS}.

In this notation, Liv\v{s}ic's remark (see \cite[Remark~1]{Liv71} and \cite[\S3, Remark~2]{Liv72}) is that, in the setting of his theorem, if $\omega$ is Dini then every $f\in\Lip_\omega(X)$ whose sums over periodic orbits vanish is a continuous coboundary; he also noted that continuity of $f$ alone is not sufficient (see \cite[\S3, Remark~2]{Liv72}).
Indeed, for a homeomorphism $T$, every continuous function whose integral vanishes for every $\mu\in\M_T$ is a continuous coboundary if and only if $T^k=\mathrm{id}$ for some $k\ge1$ (see \cite{Ko}).
In unpublished work, Conze \& Guivarc'h \cite{CG} showed that, on one-sided full shifts, every function of summable variation has a continuous sub-action.
More generally, Bousch \cite{Bo} showed that, for transitive hyperbolic systems (see Section~\ref{sec:consequences}), Liv\v{s}ic's theorem and the Ma\~n\'e lemma hold for all functions $f$ satisfying the \emph{Walters condition} \cite[D\'efinition~1]{Bo}: for all $\varepsilon>0$ there exists $\eta>0$ such that $|S_nf(x)-S_nf(y)|\le\varepsilon$ whenever $n\ge1$ and $d(T^ix,T^iy)\le\eta$ for $0\le i<n$. This condition was introduced by Walters \cite{Wa78} for maps which expand distances, in terms of their inverse branches, and he observed that it holds for functions of summable variation on one-sided subshifts of finite type \cite[\S3.3]{Wa78}; for hyperbolic systems, Bousch noted that it holds for every function in $\Lip_\omega(X)$ when $\omega$ is Dini \cite[D\'efinition-Proposition~2]{Bo}. The Walters condition implies the \emph{Bowen property}, introduced by Bowen \cite{Bow} in his work on uniqueness of equilibrium states: there are $\varepsilon>0$ and $R\ge0$ such that $|S_nf(x)-S_nf(y)|\le R$ whenever $n\ge1$ and $d(T^ix,T^iy)\le\varepsilon$ for $0\le i<n$.

Conversely, for Lipschitz uniformly expanding open maps $T$ of infinite compact metric spaces $X$, Jenkinson, Zhang \& Zhao \cite[Theorem~B]{JZZ} showed that if $\omega$ does not satisfy the Dini condition, then $\Lip_\omega(X)$ contains functions which do not have any bounded sub-action, and in particular do not have any continuous sub-action.\footnote{In \cite{JZZ}, moduli of continuity are \emph{metric transforms}, i.e.\ continuous, concave, strictly increasing functions vanishing at $0$. A concave modulus coincides with such a function on some interval $[0,a]$ with $a>0$. As $X$ is compact, the space $\Lip_\omega(X)$, the Dini condition and the behaviour of $D\omega$ near $0$ depend only on the restriction of $\omega$ to $[0,a]$, so the results of \cite{JZZ} apply to concave moduli.} These functions are constructed near a periodic point, using a criterion for maps with a \emph{weakly repelling} fixed point \cite[Theorem~C]{JZZ}; such a fixed point $x^*$ satisfies, in particular, $d(Tx,x^*)>d(x,x^*)$ for all $x\ne x^*$ in a ball centred at $x^*$. This local construction does not apply to invertible hyperbolic systems: a fixed point of an Anosov diffeomorphism, for example, is not weakly repelling, as the points of its stable manifold approach it.\footnote{The conclusion of \cite[Theorem~B]{JZZ} can nevertheless be transferred to some invertible systems, such as two-sided subshifts of finite type, 
by means of a Lipschitz factor map; this argument does not apply to transitive Anosov diffeomorphisms of connected manifolds (see Remark~\ref{rem:factor}).} 
Moreover, the necessity of the Dini condition for \emph{Liv\v{s}ic's theorem} was left open in \cite{JZZ}, even for expanding maps of the circle.

The main result of this paper is that the Dini condition is necessary for both Liv\v{s}ic's theorem and the Ma\~n\'e lemma on every infinite basic set (on a periodic orbit, both results hold for every function). 
Here a \emph{basic set} of a $C^1$ diffeomorphism $T$ of a Riemannian manifold is a non-empty compact, locally maximal, hyperbolic $T$-invariant set $\Lambda$ such that $T|_\Lambda$ is transitive, i.e.\ for all non-empty open subsets $U,V$ of $\Lambda$ there is $n\ge1$ with $T^n(U)\cap V\ne\emptyset$. 
As $T|_\Lambda$ is transitive, a basic set $\Lambda$ is finite if and only if it is a periodic orbit.
We use the Riemannian distance, both on $\Lambda$ and on the manifold, to define the spaces $\Lip_\omega$; if the manifold is not connected, we take the distance between points in different components to be a constant, large enough that this defines a metric on $\Lambda$, and on the manifold when the latter is compact (the spaces $\Lip_\omega$ do not depend on this constant).

\begin{thm}\label{thm:basic}
Let $\Lambda$ be an infinite basic set of a $C^1$ diffeomorphism $T$ of a Riemannian manifold $M$. For a concave modulus $\omega$, the following are equivalent:
\begin{itemize}
\item[(a)] $\omega$ satisfies the Dini condition;
\item[(b)] every $f\in\Lip_\omega(\Lambda)$ with $\int f\,d\mu=0$ for all $\mu\in\M_{T|_\Lambda}$ is a continuous coboundary for $T|_\Lambda$;
\item[(c)] every $f\in\Lip_\omega(\Lambda)$ has a continuous sub-action with respect to $T|_\Lambda$;
\item[(d)] every $f\in\Lip_\omega(\Lambda)$ satisfies the Walters condition for $T|_\Lambda$;
\item[(e)] every $f\in\Lip_\omega(\Lambda)$ has the Bowen property for $T|_\Lambda$;
\item[(f)] every $f\in\Lip_\omega(\Lambda)$ has a bounded sub-action with respect to $T|_\Lambda$.
\end{itemize}
If $\omega$ does not satisfy the Dini condition, then $\Lip_\omega(\Lambda)$ contains a function $f$ such that $\int f\,d\mu=0$ for every $\mu\in\M_{T|_\Lambda}$, and $\inf u=-\infty$ for every sub-action $u$ for $f$ with respect to $T|_\Lambda$; if, moreover, $M$ is compact, then $\Lip_\omega(M)$ contains a function with the same properties for $T:M\to M$.
\end{thm}

Condition (f) of Theorem~\ref{thm:basic} can also be expressed in terms of Birkhoff sums: a continuous function $f$ has a bounded sub-action if and only if $\sup_{n\ge1}\max_x\big(S_nf(x)-n\beta(f)\big)<\infty$ (see Lemma~\ref{lem:bounded}); this is condition (B) of Morris \cite{Mo07} (cf.~\cite[Remark~8.5]{JZZ}).

Theorem~\ref{thm:basic} applies, with $\Lambda=M$, to every transitive Anosov diffeomorphism $T$ of a compact manifold $M$ of positive dimension, as $M$ is then an infinite basic set. For Anosov diffeomorphisms which are not assumed to be transitive, we have the following consequence.

\begin{cor}\label{cor:anosov}
Let $T$ be an Anosov diffeomorphism of a compact manifold $M$ of positive dimension, and let $\omega$ be a concave modulus which does not satisfy the Dini condition. Then $\Lip_\omega(M)$ contains a function $f$ such that $\int f\,d\mu=0$ for every $\mu\in\M_T$, and $\inf u=-\infty$ for every sub-action $u$ for $f$.
\end{cor}

Corollary~\ref{cor:anosov} is proved in Section~\ref{sec:consequences}, by the argument used for Theorem~\ref{thm:basic}.

In Theorem~\ref{thm:basic}, the implications (a)$\Rightarrow$(d), (d)$\Rightarrow$(b) and (d)$\Rightarrow$(c) are the results of Bousch \cite{Bo} recalled above (see also Section~\ref{sec:consequences}), (d)$\Rightarrow$(e) holds because the Walters condition implies the Bowen property, (c)$\Rightarrow$(f) holds because continuous functions on the compact set $\Lambda$ are bounded, and (e)$\Rightarrow$(f) follows from a closing argument (see Lemma~\ref{lem:closing}). 
The converse implications are the main content of Theorem~\ref{thm:basic}; for these it suffices, when $\omega$ is not Dini, to find a single function $f\in\Lip_\omega(\Lambda)$ such that $\int f\,d\mu=0$ for every invariant measure $\mu$, and $\inf u=-\infty$ for every sub-action $u$ for $f$. Such a function can be written down explicitly, for any Lipschitz map with an orbit of the kind described in the following theorem.

\begin{thm}\label{thm:main}
Let $T:X\to X$ be a Lipschitz self-map of a compact metric space, and let $(z_k)_{k\in\Z}$ be a sequence in $X$ with $Tz_k=z_{k+1}$ for all $k$, such that the closures $A$ of $\{z_k:k<0\}$ and $B$ of $\{z_k:k\ge0\}$ are disjoint. Let $\omega$ be a concave modulus, let $L>1$ be a Lipschitz constant for $T$, and for $n\ge0$ let $s_n:=\big(\omega(L^{-n})-\omega(L^{-n-1})\big)/(L^{-n}-L^{-n-1})$, the slope of the chord of $\omega$ over $[L^{-n-1},L^{-n}]$. Let $d(A,B):=\min_{a\in A,\,b\in B}d(a,b)>0$, set $\psi(x):=\min\{1,d(x,A)/d(A,B)\}$, and
\begin{equation}\label{eq:f}
f_\omega:=\sum_{n=1}^\infty b_n\,(\psi\circ T^n-\psi),\qquad\text{where }b_n:=L^{-n}(s_n-s_{n-1}).
\end{equation}
Then the series \eqref{eq:f} converges uniformly, $f_\omega\in\Lip_\omega(X)$, with $|f_\omega|_\omega\le2\max\{1,d(A,B)^{-1}\}$, and
\[
\int f_\omega\,d\mu=0\qquad\text{for every }\mu\in\M_T.
\]
Moreover, the following are equivalent:
\begin{itemize}
\item[(a)] $\omega$ satisfies the Dini condition;
\item[(b)] $f_\omega=v\circ T-v$ for some $v\in C(X)$;
\item[(c)] $f_\omega$ has a sub-action which is bounded below;
\item[(d)] the Birkhoff sums $S_Nf_\omega(z_{-N})=\sum_{i=0}^{N-1}f_\omega(z_{i-N})$, $N\ge1$, are bounded above.
\end{itemize}
For every $N\ge1$,
\begin{equation}\label{eq:closed}
S_Nf_\omega(z_{-N})=\sum_{j=1}^N\omega(L^{-j})-\frac{\omega(1)-\omega(L^{-N})}{L-1},
\end{equation}
and consequently
\begin{equation}\label{eq:growth}
\frac{1}{\log L}\int_{L^{-N}}^1\frac{\omega(r)}{r}\,dr-\frac{L}{L-1}\,\omega(1)\ \le\ S_Nf_\omega(z_{-N})\ \le\ \frac{1}{\log L}\int_{L^{-N}}^1\frac{\omega(r)}{r}\,dr;
\end{equation}
in particular, if $\omega$ does not satisfy the Dini condition, then $S_Nf_\omega(z_{-N})\to+\infty$ as $N\to\infty$.
\end{thm}

If $T$ is a homeomorphism, then $z_k=T^kz_0$ for all $k$, and the hypothesis on $(z_k)$ holds if and only if the $\alpha$- and $\omega$-limit sets of $z_0$ are disjoint (see Lemma~\ref{lem:limitsets}). An infinite basic set contains a \emph{heteroclinic point}, i.e.\ a point whose $\alpha$- and $\omega$-limit sets are distinct periodic orbits, so Theorem~\ref{thm:main} applies, and gives Theorem~\ref{thm:basic} (see Section~\ref{sec:consequences}).
The hypothesis is also satisfied by some non-invertible maps, such as expanding maps of the circle, for which Theorem~\ref{thm:main} shows that the Dini condition is necessary for Liv\v{s}ic's theorem (see Remark~\ref{rem:expanding}).

The coefficients $b_n$ are non-negative, with $\sum_nb_n\le\omega(1)$, and their first moment $\sum_nnb_n$, which is the limit of $S_Nf_\omega(z_{-N})$ as $N\to\infty$, is finite if and only if $\omega$ is Dini (see Lemma~\ref{lem:coefficients}). When $\omega$ is Dini, $f_\omega$ has a transfer function in $\Lip_{D\omega}(X)$, which does not necessarily belong to $\Lip_\omega(X)$ (see Remark~\ref{rem:regularity}).
For the doubling map $x\mapsto2x\pmod1$ on $\R/\Z$, there is no such loss of regularity if and only if $D\omega\le C'\omega$ for some $C'>0$ (see Proposition~\ref{prop:doubling}); as this condition implies that every function in $\Lip_\omega$ is H\"older continuous \cite[Proposition~9.7]{JZZ}, the loss of regularity is unavoidable whenever $\Lip_\omega$ contains a function which is not H\"older continuous.

On basic sets, the growth rate in \eqref{eq:growth} is the largest possible. For every $f\in C(X)$, the sequence $\max_xS_Nf(x)-N\beta(f)$ is non-negative and $o(N)$, as $\beta(f)=\lim_{N\to\infty}\frac1N\max_xS_Nf(x)$ (see \cite[\S2]{Bo}, and cf.~\cite[Proposition~2.2]{Je}); Bousch \cite[\S2]{Bo} remarks that the speed of this convergence cannot be made more precise without further hypotheses on $f$ or on the dynamical system. The sequence is bounded if and only if $f$ has a bounded sub-action (see Lemma~\ref{lem:bounded}). On a basic set $\Lambda$, we show that it is at most $C|f|_\omega\int_{e^{-N}}^1\omega(r)/r\,dr$ for every $f\in\Lip_\omega(\Lambda)$ and $N\ge1$, where $C$ depends only on $\Lambda$ and $T$ (cf.~\cite[Lemma~4.4]{JLS} for mixing one-sided subshifts of finite type), and that, if $\Lambda$ is infinite, 
this bound is attained, up to a factor independent of $\omega$ and $N$, by the sum of $f_\omega$ and a multiple of the coboundary $\psi\circ T-\psi$ (see Proposition~\ref{prop:optimal}).

As the series \eqref{eq:f} converges uniformly, and each of its partial sums is a continuous coboundary, $f_\omega$ is a \emph{weak coboundary}, i.e.\ a uniform limit of continuous coboundaries (see \cite{Bo23,BJ}). A continuous function is a weak coboundary if and only if its integral vanishes for every $\mu\in\M_T$ (see \cite{Bo23} and the references therein, or the proof of \cite[Lemma~3]{BJ}). 

For the doubling map, a weak coboundary which is not a continuous coboundary was constructed by Zinsmeister \cite[\S3]{BJ}. It is a lacunary Fourier series $f(x)=\sum_{n\ge0}a_ne^{2\pi i2^nx}$, with $\sum_n|a_n|<\infty$ and $\sum_na_n=0$, for which the tails $\sum_{p>n}a_p$ are not square-summable; as these tails are Fourier coefficients of any integrable solution $v$ of $f=v\circ T-v$, this equation does not have any solution in $L^2$. The modulus of continuity of $f$ is not considered in \cite{BJ}. 
Bousch \cite{Bo23} showed that a generalisation of this construction gives all weak coboundaries: if $T$ is a continuous map of a compact metric space, then the weak coboundaries are exactly the sums $\sum_{n\ge1}(g_n\circ T^n-g_n)$, where $(g_n)_{n\ge1}$ is a sequence of continuous functions with $\sum_n\|g_n\|_\infty<\infty$. Zinsmeister's example is the case $g_n(x)=a_ne^{2\pi ix}$, and $f_\omega$ is the case $g_n=b_n\psi$.

Morris \cite[proof of Proposition~3]{Mo07} used a version of Zinsmeister's construction on the full shift on two symbols, to show that generic weak coboundaries have Birkhoff sums which are not bounded above. With coefficients given by the differences $\omega(2^{-n})-\omega(2^{-n-1})$, this construction gives a function in $\Lip_\omega$, for the usual metric on the shift, which has properties analogous to those of $f_\omega$ in Theorem~\ref{thm:main}, as one sees by computing its Birkhoff sums along a heteroclinic orbit as in Section~\ref{sec:proof}. The coefficients $b_n$ differ from these because the construction in \cite{Mo07} uses a locally constant function in place of $\psi$, and no such function exists when $X$ is connected
(see Remarks~\ref{rem:coefficients} and~\ref{rem:morris}).

This paper is organised as follows. In Section~\ref{sec:proof} we prove Theorem~\ref{thm:main}. In Section~\ref{sec:consequences} we collect the properties of basic sets that will be used, and prove Theorem~\ref{thm:basic} and Corollary~\ref{cor:anosov}. In Section~\ref{sec:growth} we prove Proposition~\ref{prop:optimal} on the growth of Birkhoff sums, and Section~\ref{sec:remarks} contains further remarks, together with Proposition~\ref{prop:doubling} on the regularity of transfer functions and sub-actions for the doubling map.

\section{Proof of Theorem~\ref{thm:main}}\label{sec:proof}

Throughout this section, $\omega$ is a concave modulus and $L>1$. We first study the coefficients $b_n$. For $0\le a<b$, the \emph{chord} of $\omega$ over $[a,b]$ is the affine function on $\R$ which agrees with $\omega$ at $a$ and $b$; its slope is $(\omega(b)-\omega(a))/(b-a)$. We use two elementary consequences of concavity: each chord of $\omega$ lies below the graph of $\omega$ between its endpoints and above it outside them, and $r\mapsto\omega(r)/r$ is non-increasing on $(0,\infty)$; in particular $r\,\omega(1)\le\omega(r)$ for $0\le r\le1$, $\omega(cr)\le c\,\omega(r)$ for $c\ge1$ and $r\ge0$, and $\omega$ is subadditive. For $n\ge0$ set
\begin{equation}\label{eq:slopes}
r_n:=L^{-n},\qquad s_n:=\frac{\omega(r_n)-\omega(r_{n+1})}{r_n-r_{n+1}},\qquad b_n:=r_n(s_n-s_{n-1})\quad(n\ge1).
\end{equation}
The function $r\mapsto\min\{1,L^nr\}$ increases linearly, with slope $1/r_n$, on $[0,r_n]$, and is constant on $[r_n,\infty)$. So adding $b_n\min\{1,L^nr\}$ to a function increases its slope by $s_n-s_{n-1}$ on $(0,r_n)$, and leaves it unchanged on $(r_n,\infty)$; here $s_n-s_{n-1}$ is the decrease at $r_n$ of the slope of the function $W$ in the following lemma.

\begin{lem}\label{lem:coefficients}
Let $r_n,s_n,b_n$ be as in \eqref{eq:slopes}, and let $W:[0,1]\to\R$ be the function which vanishes at $0$, agrees with $\omega$ at each $r_n$ $(n\ge0)$, and is affine on each interval $[r_{n+1},r_n]$. Then $b_n\ge0$ for all $n\ge1$,
\begin{equation}\label{eq:tails}
\sum_{n>N}b_n=\omega(r_N)-s_Nr_N=\frac{L\,\omega(r_{N+1})-\omega(r_N)}{L-1}\qquad(N\ge0),
\end{equation}
and
\begin{equation}\label{eq:W}
W(r)=s_0r+\sum_{n=1}^\infty b_n\min\{1,L^nr\}\qquad(0\le r\le1),
\end{equation}
with $W\le\omega$ on $[0,1]$. Moreover, for every $N\ge1$,
\begin{equation}\label{eq:moment}
\sum_{n=1}^\infty\min\{n,N\}\,b_n=\sum_{j=1}^N\omega(r_j)-\frac{\omega(1)-\omega(r_N)}{L-1}
\end{equation}
and
\begin{equation}\label{eq:dinibounds}
\frac{1}{\log L}\int_{r_N}^1\frac{\omega(r)}{r}\,dr-\omega(1)\ \le\ \sum_{j=1}^N\omega(r_j)\ \le\ \frac{1}{\log L}\int_{r_N}^1\frac{\omega(r)}{r}\,dr.
\end{equation}
Consequently, $\sum_{n\ge1}nb_n=\sum_{j\ge1}\omega(r_j)-\omega(1)/(L-1)$, and $\sum_{n\ge1}nb_n<\infty$ if and only if $\omega$ satisfies the Dini condition.
\end{lem}

\begin{proof}
As $\omega$ is non-decreasing, $s_n\ge0$; as $\omega$ is concave and $[r_{n+1},r_n]$ lies to the left of $[r_n,r_{n-1}]$, we have $s_n\ge s_{n-1}$, and so $b_n\ge0$. The function $W$ has slope $s_m$ on $(r_{m+1},r_m)$, and $W\le\omega$ on $[0,1]$, since on each interval $[r_{m+1},r_m]$ it coincides with the chord of $\omega$ over that interval.

For $n\ge0$ let $c_n:=\omega(r_n)-s_nr_n$, the value at $0$ of the chord of $\omega$ over $[r_{n+1},r_n]$. This chord lies above the graph of $\omega$ at $0$ and has non-negative slope, so $0\le c_n\le\omega(r_n)$, and in particular $c_n\to0$ as $n\to\infty$. As $\omega(r_{n-1})-\omega(r_n)=s_{n-1}(r_{n-1}-r_n)$, we have $c_{n-1}-c_n=r_n(s_n-s_{n-1})=b_n$ for $n\ge1$, and summing over $n>N$ gives the first equality in \eqref{eq:tails}; the second holds because $s_Nr_N=L\big(\omega(r_N)-\omega(r_{N+1})\big)/(L-1)$.

For $N\ge1$ let $W_N(r):=s_0r+\sum_{n=1}^Nb_n\min\{1,L^nr\}$. On $(r_{m+1},r_m)$ with $m<N$, the functions $\min\{1,L^nr\}$ with $n\le m$ are increasing and those with $n>m$ are constant, so the slope of $W_N$ is $s_0+\sum_{n=1}^m(s_n-s_{n-1})=s_m$; similarly, its slope on $(0,r_N)$ is $s_N$. Thus $W_N-W$ is constant on $[r_N,1]$, and since $W_N(r_N)=s_Nr_N$, we have $W_N=W-c_N$ on $[r_N,1]$. For each $r\in(0,1]$ this holds for all large $N$, so letting $N\to\infty$ gives \eqref{eq:W}.

For \eqref{eq:moment}, note that $\sum_n\min\{n,N\}\,b_n=\sum_{j=1}^N\sum_{n\ge j}b_n$, and that $\sum_{n\ge j}b_n=\big(L\,\omega(r_j)-\omega(r_{j-1})\big)/(L-1)$ by \eqref{eq:tails}; summing over $1\le j\le N$ gives \eqref{eq:moment}, as $\sum_{j=1}^N\omega(r_{j-1})=\sum_{j=1}^N\omega(r_j)+\omega(1)-\omega(r_N)$. As $\omega$ is non-decreasing, $\omega(r_j)\log L\le\int_{r_j}^{r_{j-1}}\omega(r)/r\,dr\le\omega(r_{j-1})\log L$ for $j\ge1$, and summing over $1\le j\le N$ gives
\[
\log L\sum_{j=1}^N\omega(r_j)\ \le\ \int_{r_N}^1\frac{\omega(r)}{r}\,dr\ \le\ \log L\sum_{j=1}^N\omega(r_{j-1})\ \le\ \log L\,\Big(\omega(1)+\sum_{j=1}^N\omega(r_j)\Big),
\]
which is \eqref{eq:dinibounds}. Finally, letting $N\to\infty$ in \eqref{eq:moment} gives the formula for $\sum_nnb_n$, by monotone convergence, as $\omega(r_N)\to0$; and by \eqref{eq:dinibounds}, $\sum_j\omega(r_j)<\infty$ if and only if $\omega$ satisfies the Dini condition.
\end{proof}

We can now prove Theorem~\ref{thm:main}.

\begin{proof}[Proof of Theorem~\ref{thm:main}]
Let $K:=\max\{1,d(A,B)^{-1}\}$. The function $\psi$ takes values in $[0,1]$, is $d(A,B)^{-1}$-Lipschitz, and satisfies $\psi=0$ on $A$ and $\psi=1$ on $B$; hence
\begin{equation}\label{eq:step}
\psi(z_k)=0\ \text{ for }k<0,\qquad\psi(z_k)=1\ \text{ for }k\ge0.
\end{equation}
If $r_n,s_n,b_n,W$ are as in Lemma~\ref{lem:coefficients}, then
$b_n\ge0$ and $\sum_nb_n\le\omega(1)$ (by \eqref{eq:tails} with $N=0$), so the series \eqref{eq:f} converges uniformly, $f_\omega\in C(X)$ and $\|f_\omega\|_\infty\le\omega(1)$.

Let $x,y\in X$ and $r:=d(x,y)$. If $r>1$, then $|f_\omega(x)-f_\omega(y)|\le2\omega(1)\le2\omega(r)$. If $r\le1$, then $d(T^nx,T^ny)\le L^nr$, so $|\psi(T^nx)-\psi(T^ny)|\le K\min\{1,L^nr\}$ and $|\psi(x)-\psi(y)|\le Kr$. By Lemma~\ref{lem:coefficients},
\[
|f_\omega(x)-f_\omega(y)|\le K\sum_{n=1}^\infty b_n\min\{1,L^nr\}+Kr\sum_{n=1}^\infty b_n\le KW(r)+Kr\,\omega(1)\le2K\omega(r).
\]

For $n\ge1$ set $v_n:=\sum_{j=0}^{n-1}\psi\circ T^j$, so that $0\le v_n\le n$ and $\psi\circ T^n-\psi=v_n\circ T-v_n$. Thus the partial sums of \eqref{eq:f} are continuous coboundaries, and $\int f_\omega\,d\mu=0$ for every $\mu\in\M_T$; in particular $\beta(f_\omega)=0$.

For $j\ge1$, \eqref{eq:step} gives $\psi(z_{n-j})-\psi(z_{-j})=1$ if $n\ge j$, and $0$ otherwise. Hence $f_\omega(z_{-j})=\sum_{n\ge j}b_n$, and
\begin{equation}\label{eq:sums}
S_Nf_\omega(z_{-N})=\sum_{j=1}^N\sum_{n\ge j}b_n=\sum_{n=1}^\infty\min\{n,N\}\,b_n.
\end{equation}
So \eqref{eq:closed} follows from \eqref{eq:moment}, and \eqref{eq:growth} from \eqref{eq:closed} and \eqref{eq:dinibounds}, as $0\le\omega(1)-\omega(L^{-N})\le\omega(1)$. By monotone convergence, $S_Nf_\omega(z_{-N})$ increases to $\sum_nnb_n\in[0,\infty]$ as $N\to\infty$.

We now prove (a)$\Rightarrow$(b)$\Rightarrow$(c)$\Rightarrow$(d)$\Rightarrow$(a).

(a)$\Rightarrow$(b): By Lemma~\ref{lem:coefficients}, $\sum_nnb_n<\infty$, so $v:=\sum_nb_nv_n$ converges uniformly, $v\in C(X)$, and $v\circ T-v=\sum_nb_n(\psi\circ T^n-\psi)=f_\omega$.

(b)$\Rightarrow$(c): As $\beta(f_\omega)=0$, $v$ is a sub-action, and it is bounded.

(c)$\Rightarrow$(d): If $u$ is a sub-action for $f_\omega$ which is bounded below, then summing \eqref{eq:subaction} along $z_{-N},\dots,z_{-1}$ gives $S_Nf_\omega(z_{-N})\le u(z_0)-u(z_{-N})\le u(z_0)-\inf u$.

(d)$\Rightarrow$(a): By \eqref{eq:sums}, $\sum_nnb_n<\infty$, so $\omega$ satisfies the Dini condition by Lemma~\ref{lem:coefficients}.
\end{proof}

\begin{rem}\label{rem:coefficients}
The coefficients $b_n$ are chosen so that $f_\omega\in\Lip_\omega(X)$ for every concave modulus $\omega$: in the proof of Theorem~\ref{thm:main}, the estimate $|\psi(T^nx)-\psi(T^ny)|\le K\min\{1,L^nd(x,y)\}$ is combined with the inequality $\sum_nb_n\min\{1,L^nr\}\le W(r)\le\omega(r)$ for $0\le r\le1$ from Lemma~\ref{lem:coefficients}. If the simpler coefficients $\tilde b_n:=\omega(L^{-n})-\omega(L^{-n-1})$ were chosen in place of $b_n$, the resulting function $\tilde f:=\sum_n\tilde b_n(\psi\circ T^n-\psi)$ would satisfy $S_N\tilde f(z_{-N})=\sum_{j=1}^N\omega(L^{-j})$, by the computation leading to \eqref{eq:sums} (cf.~\eqref{eq:closed}), so its Birkhoff sums along $(z_k)$ would again be bounded if and only if $\omega$ is Dini. However, the argument in the proof would only give $|\tilde f(x)-\tilde f(y)|\le2K\sum_n\tilde b_n\min\{1,L^nr\}$, where $r:=d(x,y)\le1$, and this does not show that $\tilde f\in\Lip_\omega(X)$, as the right-hand side can be much larger than $\omega(r)$: for $\omega(r)=r$, for example, it is of order $r\log(1/r)$.
\end{rem}

\begin{rem}\label{rem:morris}
In the construction of \cite[proof of Proposition~3]{Mo07}, described in Section~\ref{sec:intro}, 
the coefficients $\tilde b_n$ of Remark~\ref{rem:coefficients}, with $L=2$, give the function
$f=\sum_{n\ge1}\tilde b_n(\varphi\circ\sigma^n-\varphi)$, where $\sigma$ is the shift on the space of sequences $x=(x_i)_{i\ge0}$ with $x_i\in\{1,2\}$, with the metric $d(x,y)=2^{-\min\{i:\,x_i\ne y_i\}}$, and $\varphi(x)=1$ or $-1$ according as $x_0=1$ or $x_0=2$. This function lies in $\Lip_\omega$ because $\varphi$ is locally constant: the $n$th term of the series takes the same value at any two points at distance less than $2^{-n}$, so only the terms with $2^{-n}\le d(x,y)$ contribute to $|f(x)-f(y)|$, and their coefficients sum to at most $\omega(d(x,y))$. When $X$ is connected, this argument is not available, as a function which is $0$ on $A$ and $1$ on $B$ cannot be locally constant.
\end{rem}

We conclude this section with two general lemmas. The first is the reformulation of the hypothesis of Theorem~\ref{thm:main} for homeomorphisms. If $T$ is a homeomorphism, then $z_k=T^kz_0$ for all $k$, and the $\alpha$- and $\omega$-limit sets of $z_0$ are contained in $A$ and $B$ respectively, so they are disjoint if $A$ and $B$ are. The converse is the following lemma.

\begin{lem}\label{lem:limitsets}
Let $T:X\to X$ be a homeomorphism of a compact metric space, and let $z\in X$. If the $\alpha$- and $\omega$-limit sets of $z$ are disjoint, then so are the closures of $\{T^kz:k<0\}$ and $\{T^kz:k\ge0\}$.
\end{lem}

\begin{proof}
Let $P$ and $Q$ be the $\alpha$- and $\omega$-limit sets of $z$; they are non-empty, closed, and invariant under $T$ and $T^{-1}$. If $T^jz\in Q$ for some $j$, then the whole orbit of $z$, and hence its closure, lies in $Q$, so $P\subseteq Q$, which is impossible; similarly, $T^jz\notin P$ for all $j$. In particular, $z$ is not periodic, so the points $T^kz$, $k\in\Z$, are distinct. The closures of $\{T^kz:k<0\}$ and $\{T^kz:k\ge0\}$ are $\{T^kz:k<0\}\cup P$ and $\{T^kz:k\ge0\}\cup Q$, which are therefore disjoint.
\end{proof}

The second lemma characterises the continuous functions which have a bounded sub-action, by condition (B) of Morris \cite{Mo07} (cf.~\cite[Remark~8.5]{JZZ}); 
it justifies the reformulation of condition (f) of Theorem~\ref{thm:basic} given in Section~\ref{sec:intro}, and will be used in Sections~\ref{sec:consequences} and~\ref{sec:growth}.

\begin{lem}\label{lem:bounded}
Let $T:X\to X$ be a continuous map of a compact metric space, and let $f\in C(X)$. Then $f$ has a bounded sub-action if and only if
\begin{equation}\label{eq:B}
\sup_{n\ge1}\max_{x\in X}\big(S_nf(x)-n\beta(f)\big)<\infty.
\end{equation}
In that case, $f$ has a sub-action $u$ such that $0\le u\le\sup_{n\ge1}\max_{x\in X}\big(S_nf(x)-n\beta(f)\big)$.
\end{lem}

\begin{proof}
If $u$ is a bounded sub-action for $f$, then summing \eqref{eq:subaction} along orbits gives $S_nf-n\beta(f)\le u\circ T^n-u\le2\sup|u|$ for all $n\ge1$. Conversely, suppose that \eqref{eq:B} holds, and let $u(x):=\sup\{S_nf(y)-n\beta(f):n\ge0,\ T^ny=x\}$, where $S_0f:=0$. Then $u\ge0$, and $u$ is bounded above by the left-hand side of \eqref{eq:B}, which is non-negative (as $\int(S_nf-n\beta(f))\,d\mu=0$ for every $f$-maximizing measure $\mu$); moreover, $u\circ T\ge u+f-\beta(f)$, so $u$ is a sub-action for~$f$.
\end{proof}

\section{Basic sets}\label{sec:consequences}

Let $T:X\to X$ be a continuous map of a compact metric space. Bousch \cite{Bo} considers transitive continuous surjections $T$ with the following \emph{weak local product property} \cite[D\'efinition~4]{Bo}: for every $\varepsilon>0$ there exists $\eta>0$ such that if $(x_i)_{i\le0}$ and $(y_i)_{i\ge0}$ are sequences in $X$ with $Tx_i=x_{i+1}$ for $i<0$, $Ty_i=y_{i+1}$ for $i\ge0$, and $d(x_0,y_0)\le\eta$, then there is a sequence $(w_i)_{i\in\Z}$ with $Tw_i=w_{i+1}$, $d(w_i,x_i)\le\varepsilon$ for $i\le0$, and $d(w_i,y_i)\le\varepsilon$ for $i\ge0$. For such maps he shows that every function satisfying the Walters condition has a continuous sub-action \cite[Th\'eor\`eme~3]{Bo}, and is a continuous coboundary if its integral vanishes for every invariant measure \cite[Th\'eor\`eme~4]{Bo}.

Now let $\Lambda$ be a basic set of a $C^1$ diffeomorphism $T$. 
Then $T|_\Lambda$ is a transitive homeomorphism, and it is Lipschitz, as $T$ is $C^1$ and $\Lambda$ is compact; we fix a Lipschitz constant $L>1$ for $T|_\Lambda$. We use the following standard facts (see \cite[\S\S6.4 and 18.1]{KH}):
\begin{itemize}
\item[(H1)] there are $\delta>0$, $C_0\ge1$ and $q\in(0,1)$ such that if $x,y\in\Lambda$ and $d(T^jx,T^jy)\le\delta$ for all $j\ge0$, then $d(T^jx,T^jy)\le C_0q^j\delta$ for all $j\ge0$, and the same holds with $T^{-1}$ in place of $T$;
\item[(H2)] for every $\delta'>0$ there is $\varepsilon'>0$ such that every $\varepsilon'$-pseudo-orbit in $\Lambda$, i.e.\ every sequence $(y_j)_{j\in\Z}$ in $\Lambda$ with $d(Ty_j,y_{j+1})\le\varepsilon'$ for all $j$, is $\delta'$-shadowed by some $w\in\Lambda$, i.e.\ $d(T^jw,y_j)\le\delta'$ for all $j\in\Z$ (this uses local maximality).
\end{itemize}
By (H1), if $x,y\in\Lambda$ and $d(T^jx,T^jy)\le\delta$ for all $j\in\Z$, then $x=y$: for each $k\ge0$, (H1) with $T^{-k}x$ and $T^{-k}y$ in place of $x$ and $y$, and with $j=k$, gives $d(x,y)\le C_0q^k\delta$.
 By (H2), $T|_\Lambda$ has the weak local product property: given $\varepsilon>0$, let $\varepsilon'$ be as in (H2) for $\delta'=\varepsilon/2$, and let $\eta:=\min\{\varepsilon/2,\varepsilon'/L\}$; for $(x_i)_{i\le0}$ and $(y_i)_{i\ge0}$ as in the definition, the sequence equal to $x_i$ for $i\le0$ and to $y_i$ for $i\ge1$ is an $\varepsilon'$-pseudo-orbit, and $w_i:=T^iw$, for $w$ as in (H2), has the required properties.

 The following lemma, which will be used in the proofs of Lemma~\ref{lem:walters} and Proposition~\ref{prop:optimal}, shows that two orbit segments which stay close are exponentially close away from their ends.

\begin{lem}\label{lem:hyp}
Let $q$ be as in \textup{(H1)}. There are $\varepsilon>0$ and $C_1\ge1$ such that if $x,y\in\Lambda$, $N\ge1$ and $d(T^jx,T^jy)\le\varepsilon$ for $0\le j<N$, then $d(T^jx,T^jy)\le C_1(q^j+q^{N-1-j})$ for $0\le j<N$.
\end{lem}

\begin{proof}
Let $\delta$ and $C_0$ be as in (H1), and let $\varepsilon\in(0,\delta/2]$ be as in (H2) for $\delta'=\delta/2$. Given $x$, $y$ and $N$ as in the statement, the sequence equal to $T^jy$ for $j<N-1$ and to $T^jx$ for $j\ge N-1$ is an $\varepsilon$-pseudo-orbit in $\Lambda$; let $w\in\Lambda$ be as in (H2). Then $d(T^jw,T^jx)\le\delta/2+\varepsilon\le\delta$ for all $j\ge0$, and $d(T^jw,T^jy)\le\delta$ for all $j\le N-1$. By (H1), applied to $w$ and $x$, and to $T^{N-1}w$ and $T^{N-1}y$ with $T^{-1}$ in place of $T$, we have $d(T^jw,T^jx)\le C_0q^j\delta$ for $j\ge0$, and $d(T^jw,T^jy)\le C_0q^{N-1-j}\delta$ for $j\le N-1$. The lemma follows, with $C_1:=\max\{1,C_0\delta\}$.
\end{proof}

The following lemma is stated for hyperbolic systems in \cite[D\'efinition-Proposition~2]{Bo}; we include the short proof for completeness.

\begin{lem}\label{lem:walters}
If $\omega$ is Dini, then every $f\in\Lip_\omega(\Lambda)$ satisfies the Walters condition for $T|_\Lambda$.
\end{lem}

\begin{proof}
Let $q$ be as in (H1), and let $\varepsilon$ and $C_1$ be as in Lemma~\ref{lem:hyp}. If $0<\eta\le\varepsilon$, $x,y\in\Lambda$, $n\ge1$ and $d(T^ix,T^iy)\le\eta$ for $0\le i<n$, then $d(T^ix,T^iy)\le\min\{\eta,2C_1q^{\min\{i,n-1-i\}}\}$ for $0\le i<n$, by Lemma~\ref{lem:hyp}, so
\[
|S_nf(x)-S_nf(y)|\le2|f|_\omega\sum_{j=0}^\infty\omega\big(\min\{\eta,2C_1q^j\}\big).
\]
The right-hand side tends to $0$ as $\eta\to0$, by dominated convergence, since $\omega(2C_1q^j)\le2C_1\,\omega(q^j)$ and, as $\omega$ is non-decreasing and Dini, $\sum_{j\ge1}\omega(q^j)\le(\log(1/q))^{-1}\int_0^1\omega(r)/r\,dr<\infty$.
\end{proof}

The next lemma provides periodic points whose orbits stay close to given orbit segments, and heteroclinic points.

\begin{lem}\label{lem:orbits}
\textup{(a)} For all $\varepsilon>0$ there is $m_0\ge1$ such that for all $x\in\Lambda$ and $N\ge1$ there are $1\le m\le m_0$ and $p\in\Lambda$ with $T^{N+m}p=p$ and $d(T^jp,T^jx)\le\varepsilon$ for $0\le j<N$.

\textup{(b)} If $\Lambda$ is infinite, then $\Lambda$ contains a heteroclinic point, i.e.\ a point $z$ whose $\alpha$- and $\omega$-limit sets are distinct periodic orbits $P$ and $Q$; moreover, $z$ can be chosen so that $d(T^{-n}z,T^{-n}a)\le C_0\delta q^n$ for all $n\ge0$ and some $a\in P$, where $\delta$, $C_0$ and $q$ are as in \textup{(H1)}.
\end{lem}

\begin{proof}
Let $\delta$, $C_0$ and $q$ be as in (H1). Given $\delta'\in(0,\delta/2]$, let $\varepsilon'\in(0,\delta']$ be as in (H2). As $T|_\Lambda$ is transitive and $\Lambda$ is compact, there is $m_0\ge1$ such that for all $a,b\in\Lambda$ there are $a'\in\Lambda$ and $1\le m\le m_0$ with $d(a',a)\le\varepsilon'$ and $d(T^ma',b)\le\varepsilon'$: cover $\Lambda$ by finitely many non-empty open subsets of diameter at most $\varepsilon'$, and apply transitivity to each pair of them.

(a) Take $\delta':=\min\{\varepsilon,\delta/2\}$. Given $x\in\Lambda$ and $N\ge1$, choose $a'$ and $m$ as above for $a:=T^Nx$ and $b:=x$. The sequence of period $N+m$ which runs through $x,Tx,\dots,T^{N-1}x$ and then $a',Ta',\dots,T^{m-1}a'$ is an $\varepsilon'$-pseudo-orbit in $\Lambda$; let $p\in\Lambda$ be as in (H2). Then $T^{N+m}p$ also $\delta'$-shadows this pseudo-orbit, so $d(T^{N+m+j}p,T^jp)\le2\delta'\le\delta$ for all $j\in\Z$, and therefore $T^{N+m}p=p$. Moreover, $d(T^jp,T^jx)\le\delta'\le\varepsilon$ for $0\le j<N$.

(b) By (a), with $N=1$, periodic points are dense in $\Lambda$. As $\Lambda$ is infinite, there are periodic points $a,b\in\Lambda$ whose orbits $P$ and $Q$ are distinct. Take $\delta':=\delta/2$, and choose $a'$ and $m$ as above for these $a$ and $b$. The sequence $(y_j)_{j\in\Z}$ with $y_j:=T^ja$ for $j<0$, $y_j:=T^ja'$ for $0\le j<m$, and $y_j:=T^{j-m}b$ for $j\ge m$, is an $\varepsilon'$-pseudo-orbit in $\Lambda$; let $z\in\Lambda$ be as in (H2). Then $d(T^{-i}z,T^{-i}a)\le\delta'+\varepsilon'\le\delta$ and $d(T^{m+i}z,T^ib)\le\delta$ for all $i\ge0$, so by (H1) these distances are at most $C_0\delta q^i$. As $a$ and $b$ are periodic, it follows that the $\alpha$-limit set of $z$ is $P$, and its $\omega$-limit set is $Q$.
\end{proof}

The following consequence of Lemma~\ref{lem:orbits}(a), in which an orbit segment is shadowed by a periodic orbit, is the closing argument that will be used in the proofs of Theorem~\ref{thm:basic} and Proposition~\ref{prop:optimal}. For $f\in C(\Lambda)$ we write $\beta(f)$ for $\beta(\Lambda,T|_\Lambda,f)$.

\begin{lem}\label{lem:closing}
Let $\varepsilon>0$, and let $m_0$ be as in Lemma~\ref{lem:orbits}(a) for this $\varepsilon$. For all $f\in C(\Lambda)$, $x\in\Lambda$ and $N\ge1$, there are $1\le m\le m_0$ and $p\in\Lambda$ with $T^{N+m}p=p$ and $d(T^jp,T^jx)\le\varepsilon$ for $0\le j<N$, such that
\[
S_Nf(p)\le N\beta(f)+m_0\big(\beta(f)-\min f\big).
\]
\end{lem}

\begin{proof}
Let $m$ and $p$ be as in Lemma~\ref{lem:orbits}(a). As $T^{N+m}p=p$, the measure $\frac{1}{N+m}\sum_{j=0}^{N+m-1}\delta_{T^jp}$ is $T|_\Lambda$-invariant, so $S_{N+m}f(p)\le(N+m)\beta(f)$. Hence $S_Nf(p)=S_{N+m}f(p)-S_mf(T^Np)\le N\beta(f)+m\big(\beta(f)-\min f\big)$, and the lemma follows, as $m\le m_0$ and $\beta(f)\ge\min f$.
\end{proof}

We can now prove Theorem~\ref{thm:basic}.

\begin{proof}[Proof of Theorem~\ref{thm:basic}]
We prove (a)$\Rightarrow$(d)$\Rightarrow$(b)$\Rightarrow$(a), (d)$\Rightarrow$(c)$\Rightarrow$(f)$\Rightarrow$(a) and (d)$\Rightarrow$(e)$\Rightarrow$(f).

(a)$\Rightarrow$(d) is Lemma~\ref{lem:walters}. 

(d)$\Rightarrow$(b) and (d)$\Rightarrow$(c): As $T|_\Lambda$ is a transitive homeomorphism with the weak local product property, these follow from \cite[Th\'eor\`emes~4 and~3]{Bo} respectively. 

(d)$\Rightarrow$(e): If $f$ satisfies the Walters condition, and $\eta$ is as in that condition for $\varepsilon=1$, then $f$ has the Bowen property, with $\eta$ and $1$ in place of $\varepsilon$ and $R$.

(c)$\Rightarrow$(f) holds because continuous functions on $\Lambda$ are bounded.

(e)$\Rightarrow$(f): Suppose that $f\in\Lip_\omega(\Lambda)$ has the Bowen property, with constants $\varepsilon$ and $R$, and let $m_0$ be as in Lemma~\ref{lem:orbits}(a) for this $\varepsilon$. For $x\in\Lambda$ and $N\ge1$, let $p$ be as in Lemma~\ref{lem:closing}; then
\[
S_Nf(x)-N\beta(f)\le S_Nf(x)-S_Nf(p)+m_0\big(\beta(f)-\min f\big)\le R+m_0(\max f-\min f).
\]
So $f$ has a bounded sub-action, by Lemma~\ref{lem:bounded}.

(b)$\Rightarrow$(a) and (f)$\Rightarrow$(a): By Lemma~\ref{lem:orbits}(b), $\Lambda$ contains a heteroclinic point $z$; its $\alpha$- and $\omega$-limit sets are distinct periodic orbits, and so are disjoint. Let $z_k:=T^kz$ for $k\in\Z$. By Lemma~\ref{lem:limitsets}, applied to the homeomorphism $T|_\Lambda$, the closures $A$ of $\{z_k:k<0\}$ and $B$ of $\{z_k:k\ge0\}$ are disjoint. For each concave modulus $\omega$, let $f_\omega\in\Lip_\omega(\Lambda)$ be the function given by Theorem~\ref{thm:main} for $T|_\Lambda$, $L$ and $(z_k)_{k\in\Z}$. Then $\int f_\omega\,d\mu=0$ for every $\mu\in\M_{T|_\Lambda}$, and if $\omega$ is not Dini, then $\inf u=-\infty$ for every sub-action $u$ for $f_\omega$ with respect to $T|_\Lambda$. In that case, $f_\omega$ does not have any bounded sub-action, and it is not a continuous coboundary, as a continuous coboundary $v\circ T-v$ has the bounded sub-action $v$; so (b) and (f) do not hold.

If $M$ is compact, then $T:M\to M$ is Lipschitz, and $A$ and $B$ are also the closures of $\{z_k:k<0\}$ and $\{z_k:k\ge0\}$ in $M$, as $\Lambda$ is closed. So Theorem~\ref{thm:main}, applied to $T:M\to M$ and the same sequence, gives a function in $\Lip_\omega(M)$ with the same properties for $T:M\to M$.
\end{proof}

Corollary~\ref{cor:anosov} is proved in the same way, except that the heteroclinic point is found without using transitivity.

\begin{proof}[Proof of Corollary~\ref{cor:anosov}]
As $M$ is a compact hyperbolic set for $T$, and is trivially locally maximal, \textup{(H1)} and \textup{(H2)} hold with $\Lambda=M$, and $T$ is Lipschitz. The periodic points of $T$ are dense in its non-wandering set $\Omega$ (see \cite[Corollary~6.4.19]{KH}), and we first show that $\Omega$ is infinite. Otherwise, $\Omega$ would consist of finitely many hyperbolic periodic orbits, the $\omega$-limit set of each point of $M$ would be one of them (as an $\omega$-limit set is not the union of two disjoint non-empty closed invariant sets), and so $M$ would be the union of the stable manifolds of finitely many periodic points. This is impossible, as these stable manifolds are injectively immersed submanifolds of dimension $\dim M-\dim E^u<\dim M$, and so have measure zero. Here $\dim E^u\ge1$ on each component $M_0$ of $M$: some iterate $T^k$ maps $M_0$ onto itself, and if $E^u$ were trivial on $M_0$, then a further iterate of $T^k$ would be a strict contraction of $M_0$ onto itself, which is impossible as $M_0$ is compact and has positive diameter.

So $T$ has infinitely many periodic orbits. Let $\delta$ be as in \textup{(H1)}, and let $\varepsilon'$ be as in \textup{(H2)} for $\delta'=\delta/2$. As $M$ is compact, there are periodic points $a$ and $b$, with distinct orbits $P$ and $Q$, such that $d(a,b)\le\varepsilon'$. The sequence equal to $T^ja$ for $j<0$, and to $T^jb$ for $j\ge0$, is an $\varepsilon'$-pseudo-orbit; let $z\in M$ be as in \textup{(H2)}. As in the proof of Lemma~\ref{lem:orbits}(b), it follows from \textup{(H1)} that the $\alpha$- and $\omega$-limit sets of $z$ are $P$ and $Q$. By Lemma~\ref{lem:limitsets}, the sequence $z_k:=T^kz$ satisfies the hypothesis of Theorem~\ref{thm:main} for $T:M\to M$, and the function $f_\omega\in\Lip_\omega(M)$ given by that theorem has the required properties.
\end{proof}

\section{Growth of Birkhoff sums}\label{sec:growth}

We now show that, on basic sets, the growth rate in \eqref{eq:growth} is, up to a constant factor, the largest possible for $\max_xS_Nf(x)-N\beta(f)$ with $f\in\Lip_\omega$. For a concave modulus $\omega$ and $N\ge1$, let
\[
J_N(\omega):=\int_{e^{-N}}^1\frac{\omega(r)}{r}\,dr.
\]
As in Section~\ref{sec:consequences}, we write $\beta(f)$ for $\beta(\Lambda,T|_\Lambda,f)$. The upper bound \eqref{eq:upper} below is proved using the closing argument of Lemma~\ref{lem:closing};
for mixing one-sided subshifts of finite type, a bound of this form was proved in \cite[Lemma~4.4]{JLS}. When $\omega$ is Dini, it shows that the sums $S_Nf-N\beta(f)$ are bounded above, uniformly in $N$; this was proved by Conze \& Guivarc'h \cite{CG} for functions of summable variation on one-sided full shifts, and by Bousch \cite[Lemme~1]{Bo} for functions satisfying the Walters condition, on transitive systems with the weak local product property.

\begin{prop}\label{prop:optimal}
Let $\Lambda$ be a basic set of a $C^1$ diffeomorphism $T$, and let $L>1$ be a Lipschitz constant for $T|_\Lambda$. There is $C\ge1$ such that for every concave modulus $\omega$, every $f\in\Lip_\omega(\Lambda)$ and every $N\ge1$,
\begin{equation}\label{eq:upper}
\max_{x\in\Lambda}S_Nf(x)-N\beta(f)\le C\,|f|_\omega\,J_N(\omega).
\end{equation}
If $\Lambda$ is infinite, then $C$ can be chosen such that, for every concave modulus $\omega$ and every $N\ge1$,
\begin{equation}\label{eq:twosided}
C^{-1}J_N(\omega)\ \le\ \sup\Big\{\max_{x\in\Lambda}S_Nf(x)-N\beta(f):f\in\Lip_\omega(\Lambda),\ |f|_\omega\le1\Big\}\ \le\ C\,J_N(\omega).
\end{equation}
More precisely, if $f_\omega$ and $\psi$ are as in Theorem~\ref{thm:main}, for $T|_\Lambda$, $L$ and the orbit $(z_k)_{k\in\Z}$ of a heteroclinic point of $\Lambda$, then the function
\[
g_\omega:=f_\omega+\frac{L}{L-1}\,\omega(1)\,(\psi\circ T-\psi)
\]
belongs to $\Lip_\omega(\Lambda)$, with $|g_\omega|_\omega$ bounded independently of $\omega$, its integral vanishes for every $\mu\in\M_{T|_\Lambda}$, and for every $N\ge1$,
\begin{equation}\label{eq:lower}
S_Ng_\omega(z_{-N})=\sum_{j=0}^N\omega(L^{-j})+\frac{\omega(L^{-N})}{L-1}\ \ge\ \frac{J_N(\omega)}{\max\{1,\log L\}}.
\end{equation}
\end{prop}

\begin{proof}
For $0<\rho\le1$ and $\theta\ge1$, the substitution $r=s^\theta$ gives
\begin{equation}\label{eq:rescale}
\int_{\rho^\theta}^1\frac{\omega(r)}{r}\,dr=\theta\int_\rho^1\frac{\omega(s^\theta)}{s}\,ds\le \theta\int_\rho^1\frac{\omega(s)}{s}\,ds.
\end{equation}
Let $q$ be as in (H1), let $\varepsilon$ and $C_1$ be as in Lemma~\ref{lem:hyp}, and $m_0$ as in Lemma~\ref{lem:orbits}(a) for this $\varepsilon$. Let $x\in\Lambda$ and $N\ge1$, and let $p$ be as in Lemma~\ref{lem:closing}. Then
\[
S_Nf(p)\le N\beta(f)+m_0\big(\beta(f)-\min f\big)\le N\beta(f)+m_0|f|_\omega\,\omega(\diam\Lambda).
\]
By Lemma~\ref{lem:hyp}, as $\omega$ is subadditive and $\omega(C_1r)\le C_1\,\omega(r)$,
\[
S_Nf(x)-S_Nf(p)\le|f|_\omega\sum_{j=0}^{N-1}\omega\big(C_1q^j+C_1q^{N-1-j}\big)\le2C_1|f|_\omega\sum_{j=0}^{N-1}\omega(q^j),
\]
and as $t\mapsto\omega(q^t)$ is non-increasing, $\sum_{j=0}^{N-1}\omega(q^j)\le\omega(1)+\int_0^N\omega(q^t)\,dt=\omega(1)+c^{-1}\int_{q^N}^1\omega(r)/r\,dr$, where $c:=\log(1/q)$. The last integral is at most $\max\{1,c\}J_N(\omega)$, by \eqref{eq:rescale} with $\rho=e^{-N}$ and $\theta=c$ if $c>1$, and trivially if $c\le1$. Moreover, $J_N(\omega)\ge(1-e^{-1})\,\omega(1)\ge\omega(1)/2$, since $\omega(r)/r\ge\omega(1)$ for $0<r\le1$, and $\omega(\diam\Lambda)\le\max\{1,\diam\Lambda\}\,\omega(1)$. Combining these estimates gives \eqref{eq:upper}, with $C:=2C_1\big(2+\max\{1,c^{-1}\}\big)+2m_0\max\{1,\diam\Lambda\}$.

Now suppose that $\Lambda$ is infinite, and let $(z_k)_{k\in\Z}$, $f_\omega$, $\psi$ and $g_\omega$ be as in the statement; as in the proof of Theorem~\ref{thm:basic}, the hypothesis of Theorem~\ref{thm:main} is satisfied. As $\psi\circ T-\psi$ is a continuous coboundary, Theorem~\ref{thm:main} shows that $\int g_\omega\,d\mu=0$ for every $\mu\in\M_{T|_\Lambda}$, and so $\beta(g_\omega)=0$. Let $K:=\max\{1,d(A,B)^{-1}\}$. If $x,y\in\Lambda$ and $r:=d(x,y)$, then $|\psi(x)-\psi(y)|\le K\min\{1,r\}$ and $|\psi(Tx)-\psi(Ty)|\le K\min\{1,Lr\}$, where $\omega(1)\min\{1,t\}\le\omega(t)$ for all $t\ge0$, and $\omega(Lr)\le L\,\omega(r)$. Hence $\omega(1)\,(\psi\circ T-\psi)$ belongs to $\Lip_\omega(\Lambda)$, with seminorm at most $(L+1)K$, and as $|f_\omega|_\omega\le2K$,
\[
|g_\omega|_\omega\le K_1:=\Big(2+\frac{L(L+1)}{L-1}\Big)K.
\]
By \eqref{eq:step}, $S_N(\psi\circ T-\psi)(z_{-N})=\psi(z_0)-\psi(z_{-N})=1$, so the equality in \eqref{eq:lower} follows from \eqref{eq:closed}. By \eqref{eq:dinibounds}, $\sum_{j=0}^N\omega(L^{-j})\ge I_N/\log L$, where $I_N:=\int_{L^{-N}}^1\omega(r)/r\,dr$. If $L\ge e$, then $I_N\ge J_N(\omega)$, and otherwise \eqref{eq:rescale}, with $\rho=L^{-N}$ and $\theta=1/\log L$, gives $J_N(\omega)\le I_N/\log L$; in both cases $I_N/\log L\ge J_N(\omega)/\max\{1,\log L\}$, which proves \eqref{eq:lower}. Finally, the function $f:=g_\omega/K_1$ satisfies $|f|_\omega\le1$, $\beta(f)=0$ and $\max_{x\in\Lambda}S_Nf(x)\ge J_N(\omega)/(K_1\max\{1,\log L\})$, so \eqref{eq:twosided} holds whenever $C$ is as in \eqref{eq:upper} and $C\ge K_1\max\{1,\log L\}$.
\end{proof}

Dividing \eqref{eq:upper} by $N$ shows that $\frac1N\max_{x\in\Lambda}S_Nf(x)$ converges to $\beta(f)$ from above at rate $O(J_N(\omega)/N)$, uniformly for $|f|_\omega\le1$, and by \eqref{eq:twosided} this rate is optimal when $\Lambda$ is infinite (for the rate at which $\beta(f)$ is approximated from below by averages over periodic orbits, see \cite{BQ,CM} and the references therein). If $f\in\Lip_\omega(\Lambda)$ and $\int f\,d\mu=0$ for every $\mu\in\M_{T|_\Lambda}$, then applying \eqref{eq:upper} to $f$ and $-f$ gives $\max_{x\in\Lambda}|S_Nf(x)|\le C|f|_\omega J_N(\omega)$, and, when $\Lambda$ is infinite, \eqref{eq:lower} shows that this bound is attained by $g_\omega$, up to a constant factor, for every $N\ge1$. As $S_Nf_\omega(z_{-N})=S_Ng_\omega(z_{-N})-\frac{L}{L-1}\,\omega(1)$, it is also attained by $f_\omega$ for all large $N$ when $\omega$ is not Dini; the additional coboundary in $g_\omega$ cannot be omitted in general, since $f_\omega=0$ if $\omega$ is affine on $[0,1]$.

When $\omega$ is Dini, $J_N(\omega)\le\int_0^1\omega(r)/r\,dr$ for all $N$, so \eqref{eq:upper} and Lemma~\ref{lem:bounded} give, for every $f\in\Lip_\omega(\Lambda)$, a sub-action $u$ with $0\le u\le C|f|_\omega\int_0^1\omega(r)/r\,dr$; this is a quantitative form of the implication (a)$\Rightarrow$(f) in Theorem~\ref{thm:basic}.

\section{Further results and remarks}\label{sec:remarks}

We conclude with remarks on summable variation, on non-invertible maps and the results of \cite{JZZ}, and on the regularity of transfer functions and sub-actions; for the doubling map, Proposition~\ref{prop:doubling} characterises the moduli for which there is no loss of regularity.

\begin{rem}\label{rem:sumvar}
On a one-sided shift space of sequences $x=(x_i)_{i\ge0}$, with the metric $d(x,y)=2^{-\min\{i:\,x_i\ne y_i\}}$, a function $f$ has \emph{summable variation} if $\sum_{n\ge0}\mathrm{var}_n(f)<\infty$, where $\mathrm{var}_n(f):=\sup\{|f(x)-f(y)|:x_i=y_i\text{ for }0\le i<n\}$; this is a common hypothesis in symbolic dynamics (see e.g.\ \cite{Po,Wa75}). The minimal modulus of continuity $\omega_f(r):=\sup\{|f(x)-f(y)|:d(x,y)\le r\}$ of $f$ equals $\mathrm{var}_{n+1}(f)$ on $[2^{-n-1},2^{-n})$, so $\int_0^1\omega_f(r)/r\,dr=\log2\sum_{n\ge1}\mathrm{var}_n(f)$. If $f$ has summable variation, then $\omega_f\le\omega$ for the concave modulus $\omega(r):=r+\sum_{n\ge0}\big(\mathrm{var}_n(f)-\mathrm{var}_{n+1}(f)\big)\min\{1,2^nr\}$, which is Dini. So $f$ has summable variation if and only if $f\in\Lip_\omega$ for some concave modulus $\omega$ satisfying the Dini condition (cf.~\cite[Remark~9.2]{JZZ} for uniformly expanding maps).
\end{rem}

\begin{rem}\label{rem:expanding}
For the map $x\mapsto mx\pmod1$ on $\R/\Z$, where $m\in\Z$ and $|m|\ge2$, the hypothesis of Theorem~\ref{thm:main} holds with $z_k:=m^k/(m^2-1)\pmod1$ for all $k\in\Z$; for the doubling map ($m=2$), $z_{-n}=1/(3\cdot2^n)$ for $n\ge1$, and $z_k=1/3$ for even $k\ge0$, $z_k=2/3$ for odd $k\ge1$. Every uniformly expanding $C^1$ map of the circle is topologically conjugate to $x\mapsto mx\pmod1$, where $m$ is its degree \cite{Sh69}, and the hypothesis on $(z_k)$ is preserved by topological conjugacy, so it also holds for these maps.

For Lipschitz uniformly expanding open maps of infinite compact metric spaces, the necessity of the Dini condition for the existence of bounded sub-actions, and so for the Ma\~n\'e lemma, in $\Lip_\omega(X)$ was proved in \cite[Theorem~B]{JZZ}. When such a map has an orbit as in Theorem~\ref{thm:main}, as uniformly expanding $C^1$ maps of the circle do, Theorem~\ref{thm:main} shows that the Dini condition is also necessary for Liv\v{s}ic's theorem in $\Lip_\omega(X)$.
\end{rem}

\begin{rem}\label{rem:intermittent}
The coefficients $b_n$ depend on $T$ only through $L$, so Theorem~\ref{thm:main} does not detect the stronger integrability conditions which are necessary for the Ma\~n\'e lemma for maps with an indifferent fixed point (see \cite[Theorems~D and~F]{JZZ}, and \cite{GIR, Mo09} for earlier results).
For example, let $T$ be an expanding circle map of Manneville--Pomeau type $\alpha\in(0,1)$, such as $x\mapsto x+x^{1+\alpha}\pmod1$: this is a map of $\R/\Z$ which fixes $0$, has derivative greater than $1$ at every other point, and satisfies $T'(x)=1+\xi x^\alpha+o(x^\alpha)$ as $x\to0^+$, for some $\xi>0$ (see \cite{JZZ} for the precise definition). Then the condition $\int_0^1\omega(r)r^{-1-\alpha}\,dr<\infty$ is necessary and sufficient for the Ma\~n\'e lemma in $\Lip_\omega$, and sufficient for Liv\v{s}ic's theorem in $\Lip_\omega$ \cite[Theorems~D and~E]{JZZ}; Theorem~\ref{thm:main} does not decide whether it is also necessary for Liv\v{s}ic's theorem.
\end{rem}

\begin{rem}\label{rem:factor}
For two-sided subshifts of finite type, the necessity of the Dini condition for the Ma\~n\'e lemma can be deduced from \cite[Theorem~B]{JZZ}, by composing the functions given there for the one-sided shift with the projection $(x_i)_{i\in\Z}\mapsto(x_i)_{i\ge0}$, which is a Lipschitz factor map, and using Lemma~\ref{lem:bounded}. More generally, the same argument applies to a basic set $\Lambda$ of a $C^1$ diffeomorphism $T$ whenever $T|_\Lambda$ has a Lipschitz factor map onto a uniformly expanding open map of an infinite compact metric space, as is the case if $\Lambda$ is zero-dimensional (code $\Lambda$ by a Markov partition with pairwise disjoint rectangles, and give the one-sided shift the metric $L^{-\min\{i:\,x_i\ne y_i\}}$, where $L>1$ is a Lipschitz constant for $T|_\Lambda$).
It does not apply to transitive Anosov diffeomorphisms of connected manifolds: by uniform continuity and positive expansivity, a factor map from such a diffeomorphism onto a positively expansive map is constant on sufficiently small local stable manifolds, hence on each stable manifold, as stable manifolds are connected, and hence on the whole manifold, as they are dense.
As \cite[Theorem~B]{JZZ} concerns only the Ma\~n\'e lemma, this argument does not show that the Dini condition is necessary for Liv\v{s}ic's theorem.
\end{rem}

\begin{rem}\label{rem:regularity}
If $\omega$ is Dini, the transfer function $v=\sum_nb_nv_n$ in the proof of Theorem~\ref{thm:main} has modulus of continuity a multiple of $D\omega$, and in general it does not belong to $\Lip_\omega(X)$. Indeed, let $K$ be as in that proof. If $x,y\in X$ and $r:=d(x,y)$, then, as $t\mapsto\min\{1,L^tr\}$ is non-decreasing,
\[
|v_n(x)-v_n(y)|\le K\sum_{j=0}^{n-1}\min\{1,L^jr\}\le K\int_{-\infty}^n\min\{1,L^tr\}\,dt=\frac{K}{\log L}\int_0^r\frac{\min\{1,L^ns\}}{s}\,ds,
\]
and multiplying by $b_n$ and summing over $n$ gives $|v(x)-v(y)|\le(K/\log L)\,D\omega(r)$, since $\sum_nb_n\min\{1,L^ns\}\le\omega(s)$ for all $s\ge0$ (by Lemma~\ref{lem:coefficients} for $s\le1$, and as $\sum_nb_n\le\omega(1)$ for $s>1$). On the other hand, \eqref{eq:step} gives $v_n(z_{-N})=\max\{n-N,0\}$ for $N\ge1$, so by \eqref{eq:tails}, as in the proof of \eqref{eq:moment},
\[
v(z_{-N})=\sum_{n>N}(n-N)\,b_n=\sum_{j>N}\omega(L^{-j})-\frac{\omega(L^{-N})}{L-1}\ \ge\ \frac{D\omega(L^{-N})}{\log L}-\frac{L}{L-1}\,\omega(L^{-N}),
\]
where the inequality holds because $\log L\sum_{j\ge N}\omega(L^{-j})\ge\int_0^{L^{-N}}\omega(r)/r\,dr$, as $\omega$ is non-decreasing.

Now let $\Lambda$ be an infinite basic set, and suppose that $f_\omega$ and $v$ are defined for $T|_\Lambda$ and $z_k:=T^kz$, where $z\in\Lambda$ is a heteroclinic point 
chosen as in the second part of Lemma~\ref{lem:orbits}(b). Then $v$ vanishes on the $\alpha$-limit set $P$ of $z$, as $\psi=0$ on $P$ and $T(P)=P$, and $d(z_{-N},P)\le C_0\delta q^N$ for $N\ge1$. If, for example, $\omega(r)=(\log(1/r))^{-p}$ for all small $r>0$, with $p>1$, then $D\omega(L^{-N})=(N\log L)^{1-p}/(p-1)$ for all large $N$, whereas $\omega(L^{-N})$ and $\omega(C_0\delta q^N)$ are of order $N^{-p}$, so $v\notin\Lip_\omega(\Lambda)$. Moreover, as $\beta(f_\omega)=0$, every continuous sub-action $u$ for $f_\omega$ satisfies $(v-u)\circ T\le v-u$. This inequality is an equality at every periodic point, as the values of $v-u$ along a periodic orbit are non-increasing, and hence constant. Since periodic points are dense in $\Lambda$, by Lemma~\ref{lem:orbits}(a), it follows that $(v-u)\circ T=v-u$, and so $v-u$ is constant, as $T|_\Lambda$ is transitive. Hence $f_\omega$ does not have any sub-action in $\Lip_\omega(\Lambda)$; in particular, it is not of the form $u\circ T-u$ with $u\in\Lip_\omega(\Lambda)$. Thus, in both Liv\v{s}ic's theorem for Dini moduli \cite[Remark~1]{Liv71} and the Ma\~n\'e lemma, the solution cannot in general be chosen in $\Lip_\omega$.
\end{rem}

For $\omega(r)=r^\alpha$, with $0<\alpha\le1$, there is no such loss of regularity in Liv\v{s}ic's theorem (see \cite[Theorem~19.2.1]{KH}), consistently with $D\omega=\omega/\alpha$, and for $\omega(r)=r$ there is none in the Ma\~n\'e lemma either (see \cite{STY}). For uniformly expanding open maps, if $\omega$ is Dini, then every function in $\Lip_\omega(X)$ has a sub-action with modulus of continuity a multiple of $D\omega$ (see \cite[Theorem~A]{JZZ}, and cf.~\cite[Theorem~F]{JLS} for transitive maps); for one-sided full shifts, the sub-actions constructed in \cite{CG} have variations bounded by the tails of the series $\sum_k\mathrm{var}_k(f)$ (see Remark~\ref{rem:sumvar}). For the doubling map, the following result characterises the moduli for which there is no loss of regularity.

\begin{prop}\label{prop:doubling}
Let $T$ be the doubling map $x\mapsto2x\pmod1$ on $\R/\Z$, and let $\omega$ be a concave modulus which satisfies the Dini condition. The following are equivalent:
\begin{itemize}
\item[(a)] there is $C'>0$ such that $D\omega\le C'\omega$ on $(0,1]$;
\item[(b)] every $f\in\Lip_\omega(\R/\Z)$ has a sub-action in $\Lip_\omega(\R/\Z)$;
\item[(c)] every $f\in\Lip_\omega(\R/\Z)$ with $\int f\,d\mu=0$ for all $\mu\in\M_T$ is of the form $u\circ T-u$ with $u\in\Lip_\omega(\R/\Z)$.
\end{itemize}
\end{prop}

\begin{proof}
(a)$\Rightarrow$(b): By \cite[Theorem~A]{JZZ}, every $f\in\Lip_\omega(\R/\Z)$ has a sub-action in $\Lip_{D\omega}(\R/\Z)$, which is contained in $\Lip_\omega(\R/\Z)$ by (a).

(b)$\Rightarrow$(c): If $\int f\,d\mu=0$ for all $\mu\in\M_T$, then $\beta(f)=0$, so a sub-action $u\in\Lip_\omega(\R/\Z)$ for $f$ satisfies $u\circ T-u-f\ge0$. This continuous function has integral zero with respect to Lebesgue measure, which is $T$-invariant and has full support, so $f=u\circ T-u$.

(c)$\Rightarrow$(a): Let $(z_k)$ be the orbit in Remark~\ref{rem:expanding}, let $L=2$, and let $f_\omega$ and $v$ be as in Remark~\ref{rem:regularity}. By (c), $f_\omega=u\circ T-u$ for some $u\in\Lip_\omega(\R/\Z)$, so $(v-u)\circ T=v-u$, and $v-u$ is constant, as $T$ is transitive. As $v$ vanishes at the fixed point $0$, and $d(z_{-N},0)=2^{-N}/3$, it follows that $v(z_{-N})\le|u|_\omega\,\omega(2^{-N})$, and the lower bound for $v(z_{-N})$ in Remark~\ref{rem:regularity} gives $D\omega(2^{-N})\le(2+|u|_\omega)\log2\cdot\omega(2^{-N})$ for $N\ge1$. 
Hence (a) holds, as $D\omega$ is non-decreasing, $\omega(2r)\le2\,\omega(r)$, and $\omega\ge\omega(1/2)>0$ on $[1/2,1]$.
\end{proof}

By \cite[Proposition~9.7]{JZZ}, condition (a) of Proposition~\ref{prop:doubling} implies that $\omega(r)=O(r^\alpha)$ for some $\alpha>0$. So, for the doubling map, if $\Lip_\omega(\R/\Z)$ contains a function which is not H\"older continuous, then transfer functions and sub-actions cannot in general be chosen in $\Lip_\omega(\R/\Z)$; such Dini moduli are exactly those for which Liv\v{s}ic's remark is not a consequence of his theorem for H\"older continuous functions.

\section*{Acknowledgments}
The authors thank Zeng Lian for helpful discussions.
H.~Zhao is partially supported by the National Natural Science Foundation of China (No.~123B2006).


\begin{thebibliography}{99}

\bibitem{Boc} J. Bochi, Ergodic optimization of Birkhoff averages and Lyapunov exponents, in \textit{Proceedings of the International Congress of Mathematicians (ICM 2018)}, World Scientific, 2019, 1825--1846. 

\bibitem{Bo00} T. Bousch, Le poisson n'a pas d'ar\^etes, \textit{Ann. Inst. H. Poincar\'e Probab. Statist.} \textbf{36} (2000), 489--508.

\bibitem{Bo} T. Bousch, La condition de Walters, \textit{Ann. Sci. \'Ecole Norm. Sup. (4)} \textbf{34} (2001), 287--311.

\bibitem{Bo11} T. Bousch, Le lemme de Ma\~n\'e--Conze--Guivarc'h pour les syst\`emes amphidynamiques rectifiables, \textit{Ann. Fac. Sci. Toulouse Math. (6)} \textbf{20} (2011), 1--14.

\bibitem{Bo23} T. Bousch, Une repr\'esentation des cobords faibles d'un syst\`eme dynamique, \textit{Ann. Fac. Sci. Toulouse Math. (6)} \textbf{32} (2023), 817--821.

\bibitem{BJ} T. Bousch \& O. Jenkinson, Cohomology classes of dynamically non-negative $C^k$ functions, \textit{Invent. Math.} \textbf{148} (2002), 207--217.

\bibitem{Bow} R. Bowen, Some systems with unique equilibrium states, \textit{Math. Systems Theory} \textbf{8} (1974/75), 193--202.

\bibitem{BQ} X. Bressaud \& A. Quas, Rate of approximation of minimizing measures, \textit{Nonlinearity} \textbf{20} (2007), 845--853.

\bibitem{CM} D. Collier \& I. D. Morris, Approximating the maximum ergodic average via periodic orbits, \textit{Ergodic Theory Dynam. Systems} \textbf{28} (2008), 1081--1090.

\bibitem{CLT} G. Contreras, A. O. Lopes \& Ph. Thieullen, Lyapunov minimizing measures for expanding maps of the circle, \textit{Ergodic Theory Dynam. Systems} \textbf{21} (2001), 1379--1409.

\bibitem{CG} J.-P. Conze \& Y. Guivarc'h, Croissance des sommes ergodiques et principe variationnel, unpublished manuscript, circa 1993.

\bibitem{FJ01} A. Fan \& Y. Jiang, On Ruelle--Perron--Frobenius operators. I. Ruelle theorem, \textit{Comm. Math. Phys.} \textbf{223} (2001), 125--141.

\bibitem{GIR} E. Garibaldi \& I. Inoquio-Renteria, Dynamical obstruction to the existence of continuous sub-actions for interval maps with regularly varying property, \textit{Discrete Contin. Dyn. Syst.} \textbf{40} (2020), 2315--2333.

\bibitem{GK80} V. Guillemin \& D. Kazhdan, Some inverse spectral results for negatively curved 2-manifolds, \textit{Topology} \textbf{19} (1980), 301--312.

\bibitem{Je} O. Jenkinson, Ergodic optimization in dynamical systems, \textit{Ergodic Theory Dynam. Systems} \textbf{39} (2019), 2593--2618.

\bibitem{JZZ} O. Jenkinson, Y. Zhang \& H. Zhao, The convergence-divergence criterion for the Ma\~n\'e cohomology lemma and Liv\v{s}ic theorem, preprint, 2026, arXiv:2610.32484v2.

\bibitem{JLS} Y. Jiang, Z. Li \& Y. Sun, Typical periodic optimization and freezing phase transition: Dini potentials and beyond, preprint, 2026, arXiv:2610.31327v1.

\bibitem{Ka11} B. Kalinin, Liv\v{s}ic theorem for matrix cocycles, \textit{Ann. of Math. (2)} \textbf{173} (2011), 1025--1042.


\bibitem{KH} A. Katok \& B. Hasselblatt, \textit{Introduction to the modern theory of dynamical systems}, Encyclopedia of Mathematics and its Applications~54, Cambridge University Press, Cambridge, 1995.

\bibitem{Ko} A. Kocsard, On cohomological $C^0$-(in)stability, \textit{Bull. Braz. Math. Soc. (N.S.)} \textbf{44} (2013), 489--495.

\bibitem{Liv71} A. N. Liv\v{s}ic, Homology properties of $Y$-systems, \textit{Math. Notes} \textbf{10} (1971), 758--763; translated from \textit{Mat. Zametki} \textbf{10} (1971), 555--564.

\bibitem{Liv72} A. N. Liv\v{s}ic, Cohomology of dynamical systems, \textit{Math. USSR-Izv.} \textbf{6} (1972), 1278--1301; translated from \textit{Izv. Akad. Nauk SSSR Ser. Mat.} \textbf{36} (1972), 1296--1320.

\bibitem{LMM86} R. de la Llave, J. M. Marco \& R. Moriy\'on, Canonical perturbation theory of Anosov systems and regularity results for the Livsic cohomology equation, \textit{Ann. of Math. (2)} \textbf{123} (1986), 537--611.

\bibitem{LT} A. O. Lopes \& Ph. Thieullen, Sub-actions for Anosov diffeomorphisms, \textit{Ast\'erisque} \textbf{287} (2003), 135--146.

\bibitem{Mo07} I. D. Morris, A sufficient condition for the subordination principle in ergodic optimization, \textit{Bull. Lond. Math. Soc.} \textbf{39} (2007), 214--220.

\bibitem{Mo09} I. D. Morris, The Ma\~n\'e--Conze--Guivarc'h lemma for intermittent maps of the circle, \textit{Ergodic Theory Dynam. Systems} \textbf{29} (2009), 1603--1611.

\bibitem{Po} M. Pollicott, Rates of mixing for potentials of summable variation, \textit{Trans. Amer. Math. Soc.} \textbf{352} (2000), 843--853.

\bibitem{Sh69} M. Shub, Endomorphisms of compact differentiable manifolds, \textit{Amer. J. Math.} \textbf{91} (1969), 175--199.

\bibitem{STY} X. Su, Ph. Thieullen \& W. Yu, Lipschitz sub-actions for locally maximal hyperbolic sets of a $C^1$ map, \textit{Discrete Contin. Dyn. Syst.} \textbf{44} (2024), 656--677.


\bibitem{Wa75} P. Walters, Ruelle's operator theorem and $g$-measures, \textit{Trans. Amer. Math. Soc.} \textbf{214} (1975), 375--387.

\bibitem{Wa78} P. Walters, Invariant measures and equilibrium states for some mappings which expand distances, \textit{Trans. Amer. Math. Soc.} \textbf{236} (1978), 121--153.

\bibitem{Wi13} A. Wilkinson, The cohomological equation for partially hyperbolic diffeomorphisms, \textit{Ast\'erisque} \textbf{358} (2013), 75--165.

\end{thebibliography}
\end{document}